\documentclass[12pt]{article}

\usepackage[english]{babel}
\usepackage{amssymb,amsmath} 
\usepackage{url,hyperref}

\usepackage{fullpage}

\title{Some new results on consecutive equidivisible integers}
\author{Vladimir A. Letsko \thanks{%
Volgograd State Socio-Pedagogical University {\hspace{3mm} {\it e-mail:}  \href{mailto:val-etc@yandex.ru}{val-etc@yandex.ru}}}}
\begin{document} 
\maketitle 
\thispagestyle{empty} 

\begin{abstract} 
We have found maximum possible runs of consecutive positive integers each having exactly $k$ divisors for some fixed values of $k$. In addition, we exhibit the run of 10 consecutive positive integers each having exactly 12 divisors and two runs of 14 consecutive positive integers each having exactly 24 divisors. 
\end{abstract} 

\section{Introduction}

We let $\tau(n)$ denote the number of positive divisors of a positive integer $n$.
Following \cite{DE}, we say positive integers $m$ and $n$ to be {\it equidivisible} if $\tau(m)=\tau(n)$.
For each $k \in \mathbb Z^+$, let $D(k)$ be the set of positive integers which begin maximal runs of equidivisible numbers with exactly $k$ divisors. 
More precisely,
$$D(1) = \{1\} \quad \text{and} \quad D(k)=\{a \in \mathbb Z^+ \mid \tau(a)=k, \tau(a-1) \ne k\},\quad k>1.$$
We let $D(k,s)$ denote the set of numbers $a \in D(k)$  such that $\tau(a+i)=k$ 
for $0\le i \le s-1$ and $\tau(a+s)\ne k$. 
In other words, $D(k,s)$ is the set of positive integers $n$ that start runs of $s$ consecutive integers with exactly $k$ divisors. 
For example, 33 and 85 both belong to $D(4,3)$. 
Clearly, for a fixed $k$ and varying $s$, the sets $D(k,s)$ are pairwise disjoint and form a partition of $D(k)$.
For $n \in D(k)$ we write $L(n)=s$ if $n \in D(k,s)$.

It is easy to see that for each $k$, we have $D(k,s)=\varnothing$ for all large enough $s$. 
On the other hand, Erd\"os conjectured that for every $s$, $D(k,s)\ne\varnothing$ for some $k$, i.e., 
there exists sequences of consecutive equidivisible integers of any length (see \cite[Problem B18]{GUY}).

Let 
$$M(k)=\max\{s \in \mathbb Z^+\mid D(k,s) \ne \varnothing\}.$$
So $M(k)=1$ for all odd $k$ because two squares cannot be consecutive positive integers and $M(2)=2$ because $(2,3)$ is the only pair of consecutive primes. Another trivial example is $M(4) = 3$. Indeed, $L(33)=3$ and if $n>8$ is divisible by 4 then $\tau(n)\geq 4$.

Table~\ref{table1} gives known ranges for $M(k)$ for every even $k<30$, according to \cite{DE}.

\begin{table}[!ht]
\begin{center}
\begin{tabular}{|l||c|c|c|c|c|c|c|c|c|c|c|c|c|c|}
\hline \bf k & 2 & 4 & 6 & 8 & 10 & 12 & 14 & 16 & 18 & 20 & 22 & 24 & 26 & 28 \\
\hline \bf M(k) & 2  & 3 & 5 & 7 & 3 & 5-23 & 3 & 7 & 3-5 & 3-7 & 1-3 & 5-31 & 1-3 & 2-7 \\
\hline 
\end{tabular}
%\vspace{3mm}
\caption{Possible ranges of $M(k)$ for even $k<30$ according to \cite{DE}.}
\label{table1}
\end{center}
\end{table}

In 2006, Ali and Mishutka~\cite{BM} showed that $M(24) \ge 12$.

We have found exact values of $M(k)$ for 55 fixed values of $k$. Also we have shown  that $M(12)\ge 10$. Finally, we have found six runs of 13 and two runs of 14 consecutive equidivisible numbers each having exactly 24 divisors. 

\section{Results}

\subsection{Consecutive integers with $2p$ divisors}

First of all, let us consider the case $k=2p$, where $p>3$ is prime. D\"untsch and Eggleton~\cite{DE} proved that $M(2p) \le 3$. It seems that it is exactly equal to $3$.
Letting $n=3^{p-1}r_1$, $n+1=2^{p-1}r_2$, $n+2=q^{p-1}r_3$ or $n=q^{p-1}r_1$, $n+1=2^{p-1}r_2$, $n+2=3^{p-1}r_3$ for suitable primes $q$, $r_1$, $r_2$, $r_3$, we have found many examples of $n$ in $D(2p,3)$ by Chinese remainder theorem. 
Moreover, for all these $p$ we have found the smallest corresponding $n$.  
For instance, let $p=11$. Then $n=3^{10}\cdot 1765118938727$ starts the least triple for $q=5$. Since $53^{10}>n$ it is sufficient to look over the triples for $q \le 47$. It remains to note that for each of these $q$ we need to check only a small number of cases.
 
\begin{table}[!ht]
\begin{center}
\begin{tabular}{|p{6mm}|p{139mm}|}
\hline $p$ & $n$  \\
\hline \hline 5 & 7939375 \\
\hline 7 & 3388031791 \\
\hline 11 & 104228508212890623 \\
\hline 13 & 1489106237081787109375 \\
\hline 17 & 273062471666259918212890623 \\
\hline 19 & 804505911103256259918212890623 \\
\hline 23 & 490685203356467392256259918212890623 \\
\hline 29 & 6794675247932944436619977392256259918212890623 \\
\hline 31 & 329757106427071213106619977392256259918212890623 \\
\hline 37 & 4459248710164424946384890995893380022607743740081787109375 \\
\hline 41 & 3685099958690838758895720896109004106619977392256259918212890623 \\
\hline 43 & 1038001791494840815734697769103890995893380022607743740081787109375 \\
\hline 47 & 12229485870130123102579152313423230896109004106619977392256259918212 890623 \\
\hline \end{tabular}
\caption{Smallest elements $n$ of $D(2p,3))$, for prime $5 \le p \le 47$.}
\label{table3}
\end{center}
\end{table}
The correspondence between $p$ and $n$ for all prime $5\le p \le 47$ is given by the Table~\ref{table3}. In fact, we have found the necessary triples for every prime $p<200$. Full list one can see at the website of Mathematical Marathon~\cite{MM}.

It was observed that $q=5$ in all smallest triples we found. Of course, this property for smallest triples is not guaranteed for all $p$. 
While we have not yet observed a minimal triple with $q\ne 5$. The following example may serve as indirect evidence for their existence. For $p=5$, $q=19$, the smallest triple starts at $n=130,358,662,767$, while $p=5$, $q=43$ give even a smaller start of the smallest triple at $n=3,388,031,791$. 

\subsection{Exact values of $M(k)$ for $k$ divisible by $4$ and nondivisible by $3$}

Let $k$ is a positive integer such that $k>4$, $4 \mid k$ and $3 \nmid k$. If $n\equiv 4 \pmod 8$ then $3 \mid \tau(n)$. Hence it can not be equal to $k$. Therefore $M(k)\le 7$ and for every positive integer $n$ starting a run of 7 numbers each having $k$ divisors we have $n\equiv 5 \pmod 8$. 

Apparently, for the considered values $k$ $M(k)$  is exactly equal to 7. 

We show the method of constructing of required runs for $k=20$. Chinese remainder theorem provides existence of infinite number of positive integers $n$ such that

\begin{equation}
\begin{array}{l}
n=3 \cdot 13^4 \cdot q_1,\\ 
n+1=2 \cdot 17^4 \cdot q_2,\\
n+2=7^4 \cdot q_3,\\
n+3=3 \cdot 2^4\cdot q_4,\\
n+4=5^4 \cdot q_5,\\
n+5=2 \cdot 11^4 \cdot q_6,\\
n+6=3^4 \cdot q_7,         
\end{array}
\label{cond7} 
\end{equation}
and all prime factors of numbers $q_1, q_2, \dots, q_7$ are greater than $17$. 

For positive integer $n=76,043,484,008,534,356,379,398,200,621$ satisfying conditions~\eqref{cond7} numbers $q_1, q_2, q_4, q_6$ are prime and $\tau(q_3)=\tau(q_5)=\tau(q_7)=4$. Hence this $n$ starts a run of 7 numbers each having 20 divisors.

Of course, conditions~\eqref{cond7} are not necessary. Our choose of moduli is largely arbitrary.  

We used a similar technique to obtain $M(k)=7$ for some other values $k$. The correspondence between these $k$ and $n$ with $L(n)=M(k)=7$  is given by the Table~\ref{table4}.

\begin{table}[!ht]
\begin{center}
\begin{tabular}{|p{6mm}|p{136mm}|}
\hline $k$ & $n$ \\
\hline 
\hline 8 & 171893 \\
\hline 16 & 17476613 \\
\hline 20 & 76043484008534356379398200621 \\
\hline 28 & 452785996182923067361779632166688093890621 \\
\hline 32 & 788892193463818869 \\
\hline 40 & 1469311698340824775996499340503749 \\
\hline 44 & 549796909842455469360784994463197630002088937189442381988212890621 \\
\hline 52 & 116844969527144570418843086016323822512131294234905340708844790822
918212890621 \\
\hline 56 & 22553801100353754138758323632384069700595843749 \\
\hline 64 & 782810267531144296869 \\
\hline 80 & 3018228484495186382136305792833733749 \\
\hline 100 & 28507228700231793584389588114883537171392858414724415542297895106
76250621 \\
\hline 112 & 102896381882412847725564867238989575115313463843749 \\
\hline 128 & 5915712233391708437084350399869 \\
\hline \end{tabular}
\caption{Elements $n$ of $D(k,7))$ for which $L(n)=M(k)$.}
\label{table4}
\end{center}
\end{table}
 
\subsection{Exact value of $M(18)$}

We describe the case of $k=18$ in more detail. First of all, note that $\tau(n)=18$ implies that $n$ has one of the following forms: $p^{17}$, $p^8q$, $p^5q^2$, or $p^2q^2r$, where $p$, $q$, $r$ are primes. It immediately follows that if prime $p\mid n$ but $p^2\nmid n$ then $n=pa^2$ for some positive integer $a$ coprime to $p$.

Let $n$ starts a run of 5 consecutive numbers with 18 divisors each. D\"untsch and Eggleton~\cite{DE} showed that $n\equiv 1 \pmod 8$. 
Hence $n+1=2a^2$ and $n+3=4c$ where $a=pq$ or $a=p^4$, $c=r^2s$ or $c=r^5$ for some odd primes $p, q, r, s$. 
Now we can prove that $n+2$ is divisible by 3. 

We have 3 possibilities: $3\mid n$, $3\mid (n+1)$, or $3\mid (n+2)$.

If $3\mid n$ then $n=3b^2$ or $n=9b$. In the former case, we have $3b^2+1=2a^2$, which has no solutions modulo 3. In the latter case, we have $n+3=12p^2$. Considering $n+1=2a^2$ one can get $2a^2+2=12p^2$, which again has no solutions modulo 3. Therefore $n$ is not divisible by 3.

If $3\mid (n+1)$ then $n+1=18p^2$ for some prime $p>3$ (the case of $n+1=2\cdot 3^8$ is obviously impossible). Hence $n+4=3b^2$, implying that $6p^2=b^2-1$, which is impossible since the right hand side is divisible by 8 and $6p^2 \equiv 6 \pmod{8}$. Therefore $n+1$ is not divisible by 3 either.

Hence, we must have $3\mid (n+2)$, which we satisfy by assuming that $n+2=9b$. Since one of the five consecutive numbers must be divisible by 5, we let $n+4=25d$. These assumptions fulfill the divisibility requirements by small primes.

We further find it convenient to assume the following forms of the integers: $n=p_1^2p_2^2p_3$, $n+2=9p_4^2p_5$, $n+3=4p_6^2p_7$, $n+4=25p_8^2p_9$,
where $p_1, p_2, \dots, p_9$ are distinct primes greater than 5. 

Let $m_1 = p_1^2p_2^2$, $m_2 = 9p_4^2$, $m_3 = 4p_6^2$, and $m_4 = 25p_8^2$. 
We choose $p_1$, $p_2$, $p_4$, $p_6$, and $p_8$ be small primes such that the following system of quadratic congruences is solvable:
\begin{equation}
\begin{cases}
2x^2-1 \equiv 0 \pmod{m_1},\\
2x^2+1 \equiv 0 \pmod{m_2},\\ 
2x^2+2 \equiv 0 \pmod{m_3},\\
2x^2+3 \equiv 0 \pmod{m_4}
\end{cases}
\label{syst} 
\end{equation}

Since the moduli in system~\eqref{syst} are pairwise coprime and each congruence has 4 solutions, by Chinese remainder theorem the system solution 
consists of 256 classes of positive integers $\{x_i+jm\}_{j\geq 0}$, where $1\le i\le 256$ and $0\leq x_i<m=m_1m _2m_3m_4$. 
Each class can give us $n=2(x_i+jm)^2-1$ as soon as $n/m_1$, $(n+2)/m_2$, $(n+3)/m_3$, $(n+4)/m_4$ are prime
and $x_i+jm=qr$ or $x_i+jm=q^4$ for some primes $q$, $r$.

The smallest known $n$ starting 5 numbers with 18 divisors each was obtained for $p_1=7$, $p_2=17$, $p_4=11$, $p_6=13$, $p_8=29$, and $j=260$. It is equal to $$6,481,049,360,854,613,144,556,866,375,483,521.$$

Using computer search one can obtain many other numbers (for the aforementioned and other values of $m_1$, $m_2$, $m_3$, $m_4$) in $D(18,5)$.

\subsection{New lower bounds for $M(12)$ and $M(24)$}

We have found $n$ with $\tau(n)=12$ and $L(n)=10$, implying that $M(12) \ge 10$.

We have searched desired run looking over numbers n such that
\begin{equation}
\begin{array}{l}
n=19^2\cdot q_1;\\
n+1=2\cdot5^2\cdot q_2;\\
n+2=3\cdot17^2\cdot q_3;\\
n+3=2^2\cdot q_4;\\
n+4=23^2\cdot q_5;\\
n+5=2\cdot3^2\cdot q_6;\\
n+6=5\cdot7^2\cdot p;\\
n+7=2^5\cdot q_7;\\
n+8=3\cdot13^2\cdot q_8;\\
n+9=2\cdot11^2\cdot q_9. 
\end{array}
\label{cond10} 
\end{equation}
First of all, these conditions fulfill the divisibility requirements by 2, 3, 5, 7 and 11. We choose additional conditions to increase empiric probability of required number of divisors for psitive integers $q_i$. 

Let $a=5\cdot7^2$. Using Chinese remainder theorem we have found the least positive integer $p_0=7,623,414,751,537,859$, satisfying to the system of linear congruences:
$$\begin{cases}
ax^2-6 \equiv 0 \pmod{19^2},\\
ax^2-5 \equiv 0 \pmod{5^2},\\
ax^2-4 \equiv 0 \pmod{17^2},\\
ax^2-2 \equiv 0 \pmod{23^2},\\
ax^2-1 \equiv 0 \pmod{3^2},\\
ax^2+1 \equiv 0 \pmod{2^5},\\
ax^2+2 \equiv 0 \pmod{13^2},\\
ax^2+3 \equiv 0 \pmod{11^2}
\end{cases}$$

We have looked over positive integers $p=p_0+jm$, where $m=2^6\cdot3^3\cdot5\cdot11^2\cdot13^2\cdot17^2\cdot19^2\cdot23^2$. For each $p$ which is prime we can obtain $n=ap-6$ and numbers $q_1, q_2, \dots, q_9$ from conditions~\eqref{cond10}.   
It is enough to us to find $j$ for which $q_2, q_3, q_6, q_7, q_8, q_9$ are prime and $\tau(q_1)=\tau(q_4)=\tau(q_5)=4$. 
The smallest $j$ satisfying these conditions is $j=647,045,875$, which corresponds to $n=1,545,780,028,345,667,311,380,575,449$ starting a run of 10 consecutive equidivisible numbers each having exactly 12 divisors.

Using similar construction, we have also found six runs of 13 consecutive numbers each having exactly 24 divisors. Smallest of them starts at 
$$58,032,555,961,853,414,629,544,105,797,569.$$ This number gives the upper bound for \texttt{A006558$(13)$} in the OEIS~\cite{OEIS} (\texttt{A006558$(n)$} starts the first run of $n$ consecutive integers with same number of divisors).

Finally, we have found two runs of 14 consecutive equidivisible numbers each having exactly 24 divisors. 
The first run starts at $$25,335,305,376,270,095,455,498,383,578,391,968,$$ 
which was obtained in the form $3\cdot11\cdot23^2(p_0+jm)-13$ with 
$m=331,805,549,004,454,324,800$, $p_0=92,513,784,488,630,385,533$, and $j=4,373,940,659$. 
The second run starts at $$54,546,232,085,777,926,508,945,202,650,399,569,568,$$ 
which was obtained for the same values of $m$ and $p_0$ but $j=9,416,976,775,575$. Thus $M(24)\geq 14$.
Of course, the smallest of these numbers gives the upper bound for \texttt{A006558$(14)$}. 

Our results confirm the Erd\"os conjecture for $s\leq 14$.

\subsection{Smallest elements of $D(k,s)$}

\begin{table}[!ht]
\begin{center}
\begin{tabular}{|l|p{90mm}|c|c|}
\hline $k$ & $n$ & $L(n)$ & $M(k)$ \\
\hline \hline 2 & {\bf 2}  & 2 & 2 \\
\hline 4 & {\bf 33}  & 3 & 3 \\
\hline 6 &  {\bf 10093613546512321} & 5 & 5 \\
\hline 8 &  {\bf 171893} & 7 & 7 \\
\hline 10 & {\bf 7939375} & 3 & 3 \\
\hline 12 & 1545780028345667311380575449$^{\star}$ & 10 & $\le23$ \\
\hline 14 & {\bf 3388031791} & 3 & 3 \\
\hline 16 & {\bf 17476613} & 7 & 7 \\
\hline 18 & 6481049360854613144556866375483521$^{\star}$ & 5 & 5 \\
\hline 20 & 76043484008534356379398200621$^{\star}$ & 7 & 7 \\
\hline 22 & {\bf 104228508212890623}$^{\star}$ & 3 & 3 \\
\hline 24 & 25335305376270095455498383578391968$^{\star}$ & 14 & $\le31$ \\
\hline 26 & {\bf 1489106237081787109375}$^{\star}$ & 3 & 3 \\
\hline 28 & 452785996182923067361779632166688093890621$^{\star}$ & 7 & 7 \\
\hline 
\end{tabular}
\caption{Smallest known elements $n\in D(k)$ with largest $L(n)$ for even $k<30$.}
\label{table2}
\end{center}
\end{table}

Table~\ref{table2} gives a smallest known element $n\in D(k,s)$ for every even $k$ below $30$, where $s$ is greatest known number with $D(k,s) \ne \varnothing$.

Numbers $n$ which are guaranteed to be the smallest in the corresponding $D(k,L(n))$ are in bold. 
Sign ``$\star$'' in the  \texttt{$n$} column showes that corresponding $n$ have been discovered by the author of this paper.

{\bf Acknowledgement.}
We would like to thank Max A. Alekseyev who carefully read our original manuscript and proposed several corrections.

\end{document}